\documentclass[12pt]{article}
\usepackage{geometry}                
\usepackage{graphicx,color,mathtools}
\usepackage{amssymb,amsmath,amsthm,mathrsfs}
\usepackage[all,cmtip]{xy}
\usepackage{epstopdf, comment, url}
\usepackage{bm} 
\usepackage{enumitem}

\usepackage[pdftex,bookmarks,pdfnewwindow,plainpages=false,unicode,pdfencoding=auto]{hyperref}

 \hypersetup{
pdfauthor={Oleg Ivrii, Truong Vu},
pdftitle={Intersections of fractal sets with rectifiable curves},
pdfsubject={Geometric measure theory},
bookmarksdepth={4}
}

\numberwithin{equation}{section}

\newtheorem{theorem}{Theorem}[section]

\newtheorem{lemma}[theorem]{Lemma}
\newtheorem{corollary}[theorem]{Corollary}

\theoremstyle{remark}
\newtheorem*{remark}{Remark}

\theoremstyle{definition}

\DeclareMathOperator{\Hdim}{H.dim }
\DeclareMathOperator{\Adim}{A.dim }
\DeclareMathOperator{\supp}{supp }
\DeclareMathOperator{\diam}{diam}

\DeclareMathOperator{\dist}{dist}

\title{Intersections of fractal sets with rectifiable curves}
\author{Oleg Ivrii, Truong Vu}
\date{October 30, 2026}

\begin{document}

\maketitle

\begin{abstract}
Let $K \subset\mathbb{R}^d$, $d\ge2$, be a connected compact set and $\phi(t)=tL(t)$ be a gauge function with $L$ non-increasing. We study intersections of $K$ with rectifiable curves and chord-arc curves. We show that if $\int_0^1\frac{dt}{\phi(t)}<\infty$, then there exists a rectifiable curve $\gamma\subset\mathbb{R}^d$ such that $\mathcal H^\phi(K \cap\gamma)>0$. Moreover, for every $0<s<1$, there exists a chord-arc curve $\gamma$ such that $\Hdim(K \cap\gamma)\ge s$. Conversely, we show that if $K$ is uniformly porous on lines and $\int_0^1\frac{dt}{\phi(t)}=\infty$, then $\mathcal H^\phi(K \cap\gamma)=0$ for every rectifiable curve $\gamma$ and $\Hdim(K \cap\gamma)<1$ for every chord-arc curve $\gamma$.
\end{abstract}

\section{Introduction}

Although the von Koch snowflake $S$ has Hausdorff dimension
$$
\frac{\log 4}{\log 3} > 1,
$$
the intersection of $S$ with any rectifiable curve $\gamma$ has zero length. In other words, $S$ is purely unrectifiable. On the other hand, Balka and Harangi \cite{balka-harangi} showed that for any non-degenerate continuum $K$, one can find a rectifiable curve which intersects $K$ in a set of Hausdorff dimension 1.

This raises a natural question: how large can the intersection of $S$ with a rectifiable curve be? Since Hausdorff dimension alone does not distinguish between the various possibilities, we study this question in terms of Hausdorff measures with finer gauge functions. We consider gauge functions of the form $\phi(t) = t L(t) $ where $L$ is non-increasing. For the von Koch snowflake, the answer is governed by the convergence or divergence of
$$
\int_0^1 \frac{dt}{\phi(t)}.
$$
More precisely, we show that if the integral diverges, then
$$
\mathcal H^\phi (S \cap \gamma) = 0
$$
for every rectifiable curve $\gamma \subset \mathbb{C}$, whereas if the integral converges, then
there exists a rectifiable curve $\gamma$ such that
$$
\mathcal H^\phi (S \cap \gamma) > 0.
$$

In fact, the above dichotomy holds in much greater generality. 
Following Cohen \cite{cohen}, we say that a compact set $K \subset \mathbb{R}^d$ is {\em uniformly porous on lines} or UPL if there exists an $\eta > 0$ such that every non-degenerate line segment $I \subset \mathbb{R}^d$ contains a point $x \in I$ with $B(x, \eta |I|) \cap K = \varnothing$, where $|I|$ denotes the length of $I$.

The UPL condition has a natural interpretation in terms of weak tangents.
We say that $Z \subset \mathbb{R}^d$ is a {\em weak tangent} of $K$ if there are points $x_j \in K$ and scales $r_j \to 0$ such that
$$
\frac{K - x_j}{r_j} \to Z
$$
locally in the Hausdorff topology. It is easy to see that a set $K$ is UPL if no weak tangent of $K$ contains a line segment.

It is well known that a measurable set contained in a rectifiable curve has an approximate tangent at $\mathcal H^1$ almost every point. Applying this fact to $K \cap \gamma$  shows that every UPL set $K$ is purely unrectifiable. For a geometric function theory approach to pure unrectifiability, we refer the reader to \cite{joshi, younsi}.

We show:

\begin{theorem}
\label{sharp-gauges}
Let $K$ be a compact UPL set in $\mathbb{R}^d$, $d \ge 2$ and $\phi(t) = t L(t) $ be a gauge function with $L$ non-increasing. 
If $$
\int_0^1 \frac{dt}{\phi(t)} = \infty,
$$
then
$
\mathcal H^\phi (K \cap \gamma) = 0
$
for every rectifiable curve $\gamma \subset \mathbb{R}^d$. Conversely, for any non-degenerate connected compact set $K \subset \mathbb{R}^d$ and any gauge function satisfying
$$
\int_0^1 \frac{dt}{\phi(t)} < \infty,
$$
there is a rectifiable curve $\gamma$ such that
$
\mathcal H^\phi (K \cap \gamma) > 0.
$
\end{theorem}

We say that a simple curve $\gamma \subset \mathbb{R}^d$ is a {\em chord-arc curve} if for every pair of points $a, b \in \gamma$, the length of the subarc of $\gamma$ joining $a$ and $b$ is at most $C|a-b|$.
Since chord-arc curves are rather special rectifiable curves, it is natural to expect that intersections of fractal sets with chord-arc curves are smaller. The following theorem shows that this is indeed the case:

\begin{theorem}
\label{cac-sharp-gauges}
If $K$ is a compact UPL set in $\mathbb{R}^d$, $d \ge 2$, then
$$
\Hdim (K \cap \gamma) < 1
$$
for every chord-arc curve $\gamma \subset \mathbb{R}^d$. Conversely, for any non-degenerate connected compact set $K \subset \mathbb{R}^d$ and $0 < s < 1$, there is a chord-arc curve $\gamma$ such that
$$
\Hdim (K \cap \gamma) > s.
$$
\end{theorem}

\subsection*{Example 1: Self-similar sets}

To see that the von Koch snowflake $S$ is uniformly porous on lines, suppose to the contrary that $S$ contains increasingly small pieces whose rescaled copies converge to a line segment. The self-similarity of the von Koch snowflake allows us to transfer these pieces to a fixed macroscopic scale. After passing to a subsequence, we obtain a sequence of subsets of $S$ which converges in the Hausdorff topology to a line segment. As $S$ is compact, the limiting segment must be contained in $S$. This is impossible since the von Koch snowflake does not contain any line segments. 

A similar argument using the conformal elevator principle in place of exact self-similarity shows that a hyperbolic Julia set is either UPL or contains a real-analytic curve. 

\subsection*{Example 2: Uniformly wiggly quasicircles}

Let $\Sigma \subset \mathbb{R}^d$, $d \ge 2$, be a Jordan curve. We say that $\Sigma$ is a {\em quasicircle} if any two points $a, b \in \Sigma$ are contained in an arc $\sigma \subset \Sigma$ with $\diam \sigma \le C|a-b|$. 

For $x \in \Sigma$ and $0 < r < \diam \Sigma$, the {\em beta number} $\beta_\Sigma(x, r)$ measures how close $\Sigma$ is to a line segment in $B(x, r)$. More precisely,
\begin{equation}
\beta_\Sigma(x,r) := \frac{1}{r} \, \inf_L \sup_{z\in \Sigma \cap B(x,r)} \dist(z,L), 
\end{equation} 
where the infimum is taken over all lines in $\mathbb{R}^d$. We say that $\Sigma$ is {\em uniformly wiggly} if there exists a constant
$\beta_0 > 0$ such that $\beta_\Sigma(x, r) > \beta_0$ for any $x \in \Sigma$ and $0 < r < \diam \Sigma$.  By \cite[Theorem 10.2.5]{bishop-peres}, uniformly wiggly curves have Hausdorff dimension strictly greater than 1. We leave it as an exercise to the reader to check that any uniformly wiggly quasicircle is uniformly porous on lines.

\section{Background on the traveling salesman problem}

In this section, we recall sufficient conditions for a bounded set in $\mathbb R^d$ to be contained in a rectifiable or a chord-arc curve.

For a bounded set $F \subset \mathbb{R}^d$, let $\mathcal N(r)$ denote the smallest number of balls of radius $r$ needed to cover $F$. The following lemma is a simple consequence of the traveling salesman theorem of P.~Jones \cite[Theorem 10.5.1]{bishop-peres}:

\begin{lemma}
\label{jones-easy}
Let $F$ be a bounded subset of $\mathbb{R}^d$. If
$$
\int_0^{\diam F} \mathcal N(r) dr < \infty,
$$
then $F$ is contained in a rectifiable curve.
\end{lemma}

For completeness, we give an elementary direct proof:

\begin{proof}
Let $r_n = 2^{-n} \diam F$. For each $n = 0, 1, 2, \dots$, choose a finite $2r_n$-net $E_n \subset F$ with $\# E_n \le \mathcal N(r_n)$. The set $E_0$ consists of a single point, which we denote $x_0$. For every $x \in E_{n+1}$, choose a point $p(x) \in E_n$ with
$$
|x - p(x)| \le 2 r_n.
$$
Connect each $x \in E_{n+1}$ to its parent $p(x)$ by a line segment. This gives a geometric tree $\mathcal T$ whose vertices are the points in $\bigcup_n E_n$. The total length of the edges of $\mathcal T$ is bounded by
$$
\sum_{n=0}^\infty 2 r_n \, \mathcal N(r_{n+1}) \, \lesssim \, \sum_{n=0}^\infty  \int_{r_{n+2}}^{r_{n+1}} \mathcal N(r) dr \, \le \, \int_0^{\diam F} \mathcal N(r) dr.
$$
For each $n = 1, 2, \dots$, let $\mathcal T^{(n)}\subset\mathcal T$ be the finite subtree consisting of the vertices in $E_0,\dots,E_n$, and $\gamma_n$ be the unit speed curve which starts at $x_0 \in E_0$ and traverses $\mathcal T^{(n)}$ using depth-first search. Since each edge is traversed at most twice, the lengths of the curves $\gamma_n$ are uniformly bounded. The desired curve $\gamma$ is then obtained with help of the Arzelà-Ascoli theorem. Since $\gamma$ is a closed set and each point of $F$ is arbitrarily close to points in $E_n$, it follows that $\gamma$ contains $F$.
\end{proof}

A result of Badger and Vellis provides a sufficient condition for a bounded set to be contained in a chord-arc curve. The {\em Assouad dimension}\/ of a bounded set $F \subset \mathbb{R}^d$ is defined as the infimum of all $s \ge 0$ for which
$$
\mathcal N_{F \cap B(x,R)}(r)
\lesssim
 (R/r )^s,
\qquad x \in F, \qquad 0 < r < R < \diam F.
$$
It is well-known that the Hausdorff dimension of a set is at most its Assouad dimension, although the inequality may be strict.

\begin{lemma}[\textnormal{\cite[Corollary 3.5]{BV19}}] 
\label{badger-vellis}
Every bounded subset of $\mathbb R^d$ with Assouad dimension strictly smaller than $1$ is contained in a chord-arc curve.
\end{lemma}

\section{Connected compact sets}
\label{sec:compact-connected-sets}

In this section, we study intersections of compact connected sets with rectifiable and chord-arc curves.
We start with rectifiable curves:

\begin{theorem}
\label{continuum-frostman}
Let $K \subset\mathbb R^d$ be a connected compact set which contains at least two points. If
$$
\int_0^1 \frac{dt}{\phi(t)} < \infty,
$$
then there exists a rectifiable curve $\gamma \subset \mathbb{R}^d$ such that $\mathcal H^\phi(K \cap \gamma)>0$.
\end{theorem}

Our argument is a slight refinement of the one used by Balka and Harangi to show that every non-degenerate connected compact set meets some rectifiable curve in a set of Hausdorff dimension $1$. The construction is based on the following lemma:

\begin{lemma}[\textnormal{\cite[Lemma~3.3]{balka-harangi}}]
\label{cube-selection}
Suppose that $K\subset\mathbb R^d$ is a connected compact set contained in a cube $Q$ of side length 1 which meets two opposite faces of $Q$. For every integer $N \ge 2$, there are $N$ subcubes $Q_1,\ldots,Q_N\subset Q$ of side length $1/N$ with disjoint interiors such that for every $i = 1, 2, \dots, N$, there is a connected compact set $K_i\subset K\cap Q_i$ which meets two opposite faces of $Q_i$.
\end{lemma}

We will use the elementary observation
\begin{equation}
\label{eq:discretizing-the-integral}
\int_0^1 \frac{dt}{tL(t)} = \infty \quad \Longleftrightarrow \qquad \sum_{n=1}^\infty \frac{1}{L(C q^n)} = \infty,
\end{equation}
for any $C > 0$ and $0 < q < 1$. 

\begin{proof}[Proof of Theorem \ref{continuum-frostman}]
Inspection shows that there is a cube $Q_{\varnothing}$ of side length $\diam K$ which contains $K$ and whose two opposite faces meet $K$. We recursively apply Lemma \ref{cube-selection} with $N=2$. Thus, for every binary word $\omega \in\{0,1\}^n$, we obtain a cube $Q_\omega$ of side length $s_n := 2^{-n}\diam K$ and a connected compact set $K_\omega \subset K\cap Q_\omega$ meeting two opposite faces of $Q_\omega$. At each step,
$$
Q_{\omega j}\subset Q_\omega,
\qquad
K_{\omega j}\subset K_\omega \cap Q_{\omega j},
\qquad j = 0,1,
$$
and the cubes in each generation have pairwise disjoint interiors.

We define a positive measure $\mu$ on $\{0,1\}^{\mathbb{N}}$ of mass $\phi(s_0)$ by specifying it on the cylinder sets $[\omega]$ in a consistent manner. We set $\mu([\varnothing]) = \phi(s_0)$ and recursively distribute the mass among the children so that
$$
\mu([\omega]) = \mu([\omega0]) + \mu([\omega1]), \qquad \mu([\omega])  \le  \phi(s_n),
$$
for every word $\omega$ of length $n$. Since $L$ is non-increasing,
$$
\phi(s_n) \, = \, s_nL(s_n) \, \le \, 2s_{n+1}L(s_{n+1}) \, = \, 2\phi(s_{n+1}).
$$
Suppose that $\mu([\omega])$ has been defined for $\omega \in \{ 0, 1 \}^n$ and satisfies $\mu([\omega]) \le \phi(s_n)$. We distribute the mass of $[\omega]$ among its two children so that 
$$
\mu([\omega j]) \le \phi(s_{n+1}), \qquad j=0,1,
$$
and at most one child has mass strictly between 0 and $\phi(s_{n+1})$. We call such a child {\em positive non-saturated}\/. Children of mass $\phi(s_{n+1})$ are called {\em saturated}\/.

Since the side lengths $s_n \to 0$, for every infinite word  $\zeta = \zeta_1 \zeta_2 \dots \in \{ 0, 1\}^{\mathbb{N}}$, the intersection
$$
\bigcap_{n \ge 1} Q_{\zeta_1 \zeta_2 \dots \zeta_n}
$$
consists of a single point, which we denote by $\pi(\zeta)$. Define $\lambda = \pi_* \mu$ as the pushforward of $\mu$ under the coding map and let  $F:=\overline{\supp}\ \lambda$ be its closed support, that is, the smallest closed subset of $\mathbb{R}^d$ with 
$\lambda(\mathbb{R}^d \setminus F) = 0$. Since every cube $Q_{\zeta_1 \zeta_2 \dots \zeta_n}$ contains a point of $K$ and $K$ is a closed set, the image of $\pi$ is contained in $K$. Consequently, $F \subset K$.

Let $\mathcal A_n = \bigl \{ \omega \in\{0,1\}^n: \mu([\omega]) >0 \bigr \}$ be the set of words indexing the ``active'' cubes of generation $n$ which contain a point of $F$. Define $A_n = s_n\,\#\mathcal A_n$. Since each active cube of generation $n$ has at most one positive non-saturated child, the non-saturated children contribute at most $A_n/2$ to $A_{n+1}$.  On the other hand, since a saturated child $\xi$ of generation $n+1$ contributes
$$
s_{n+1}=\frac{\mu([\xi])}{L(s_{n+1})},
$$
the total contribution of the saturated children is at most
$\mu([\varnothing])/L(s_{n+1})$.  Therefore,
$$
A_{n+1} \le \frac{1}{2} \cdot A_n + \frac{\mu([\varnothing])}{L(s_{n+1})}, \qquad n = 0, 1, 2, \dots
$$
Rearranging, we get
\begin{equation}
\label{eq:continuum-forcing}
A_{n+1} - \frac{1}{2} \cdot A_n \le \frac{\mu([\varnothing])}{L(s_{n+1})}, \qquad n = 0, 1, 2, \dots
\end{equation}
By (\ref{eq:discretizing-the-integral}), the hypothesis on $\phi$ implies that
$$
\sum_{n=0}^\infty \frac{\mu([\varnothing])}{L(s_n)}<\infty.
$$
Summing (\ref{eq:continuum-forcing}) over $n = 0, 1, \dots, N$ and taking $N \to \infty$, we obtain
$$
\int_0^{\diam F} \mathcal N(r) dr \, \lesssim \, \sum_{n=0}^\infty A_n \, < \, \infty.
$$
Lemma \ref{jones-easy} produces a rectifiable curve $\gamma$ containing $F$.

By construction, $\lambda(Q_\omega) \le \phi(s_n)$ for any cube $Q_\omega$ of generation $n$.  Given $0 < r < \diam K$, choose $n$ so that
$$
s_{n+1} \, \le \, r \, < \, s_{n} \, = \, 2s_{n+1}.
$$
Since cubes of generation $n$ have side length $s_n$ and pairwise disjoint interiors, a ball of radius $r$ intersects at most $C_d$ of them.  Hence,
$$
\lambda(B(x,r)) \, \le \, C_d\,\phi(s_n) \, \le \, 2 C_d\,\phi(r).
$$
By the mass distribution principle \cite[Lemma 1.2.8]{bishop-peres}, $\mathcal H^\phi(F)>0$ as desired.
\end{proof}

\begin{remark}
Balka and Harangi \cite{balka-harangi} applied the cube selection lemma with $N$ large but retained only $N-1$ cubes at each step to produce a Cantor set in $K$ of Hausdorff dimension arbitrarily close to $1$ contained in a rectifiable curve. In the proof above, we used Lemma \ref{cube-selection} with $N=2$ but retained more of the selected cubes.
\end{remark}

\begin{theorem}
\label{continuum-frostman2}
Let $K \subset\mathbb R^d$ be a connected compact set which contains at least two points.
For every $0 < s  < 1$, there exists a chord-arc curve $\gamma \subset \mathbb{R}^d$ such that
$\Hdim (K \cap \gamma) > s.$
\end{theorem}

\begin{proof}
We follow the original construction of Balka and Harangi, where we take $N$ large and retain $N-1$ cubes at each step. Thus, for each finite word $\omega \in \{0, \dots, N-2\}^*$, we obtain a cube $Q_\omega$ of side length $N^{-|\omega|} \diam K$, where $|\omega|$ denotes the length of $\omega$. We denote the resulting limit set by
$$
F = \bigcap_{n=0}^\infty \bigcup_{|\omega| = n} Q_\omega.
$$
Since each cube $Q_\omega$ contains a point of $K$ and $K$ is a closed set, $F \subset K$. We claim that $\Hdim F = \Adim F = \log_N (N-1)$.

The lower bound for $\Hdim F$ follows after applying the mass distribution principle to the probability measure that assigns equal mass to all cubes of the same generation. Since $\Hdim F\leq\Adim F$, it remains only to prove the upper bound for the Assouad dimension.

To this end, fix a point $x \in K$ and radii $0 < r < R < \diam K$. Choose integers $n \ge m \ge 0$ so that 
\begin{align*}
N^{-m-1} \diam K \, & < \, R \, \le \, N^{-m} \diam K, \\
N^{-n-1} \diam K \, & < \, r \, \le \, N^{-n} \diam K.
\end{align*}
The ball $B(x, R)$ intersects a uniformly bounded number of cubes of generation $m$. Each such cube contains at most $(N-1)^{n-m}$ active cubes of generation $n$. Since each cube of generation $n$ can be covered by a uniformly bounded number of balls of radius $r$, we have
$$
\mathcal N_{F \cap B(x,R)}(r) \, \lesssim \, (N-1)^{n-m} \, \asymp \, \biggl ( \frac{R}{r} \biggr )^{\log_N (N-1)},
$$
as desired.  By Lemma \ref{badger-vellis}, $F$ is contained in a chord-arc curve $\gamma$. By taking $N$ sufficiently large, we can arrange that $\Hdim F > s$. Consequently,
$$
\Hdim(K \cap \gamma) \, \ge \, \Hdim F \, > \, s.
$$
The proof is complete.
\end{proof}

\section{Intersections with rectifiable curves}

For $n \ge 0$, let $\mathcal D_n$ be the collection of dyadic cubes in $[0,1]^d$ of sidelength $2^{-n}$ and $\mathcal D = \bigcup_{n=0}^\infty \mathcal D_n$.

Given a closed subset $F \subset [0,1]^d$, we say that a cube $Q\in \mathcal D$ is {\em active} if $Q \cap F \ne \varnothing$. Otherwise, we say that $Q$ is {\em passive}\/. We define the {\em dyadic energy} of $F$ as
\begin{equation}
\label{eq:energy-def}
\mathscr E(F)  =  \sum_{Q \in \mathcal D, \, Q \cap F \ne \varnothing} \diam Q.
\end{equation}

In this section, we characterize the closed subsets of a UPL set which are contained in a rectifiable curve in terms of energy and covering numbers:

\begin{theorem}
\label{energy-comparison}
Let $K$ be a compact UPL set in $[0,1]^d$, $d \ge 2$.
For a closed set $F \subset K$, the following conditions are equivalent:

{\em (1)} $F$ is contained in a rectifiable curve.

{\em (2)} $\mathscr E(F) < \infty$.

{\em (3)} $\int_0^1 \mathcal N(r) dr < \infty$.
\end{theorem}

\begin{corollary}[Divergence implies zero measure] 
\label{divergence-zero-measure}
If 
$$
\int_0^1 \frac{dt}{\phi(t)} = \infty,
$$
then $\mathcal H^\phi (K \cap \gamma) = 0$ for every rectifiable curve $\gamma$.
\end{corollary}

\begin{proof}
Set $F = K \cap \gamma$. By Theorem \ref{energy-comparison},
$$
\int_{0}^{1}  \mathcal N(r) dr < \infty.
$$
Combined with the assumption on $\phi$, this implies that
$$
\liminf_{r \to 0} \mathcal N(r) \phi(r) = 0
$$
and hence $\mathcal H^\phi(F) = 0$. 
\end{proof}

Together with Theorem \ref{continuum-frostman}, this proves Theorem \ref{sharp-gauges}. Given a line $L \subset \mathbb{R}^d$ and $r > 0$, we denote the tube of radius $r$ around $L$ by
$$
T(L, r) :=  \{ x : \dist(x, L) \le r \}.
$$

\begin{lemma}
\label{UPL-conditions}
Let $K \subset \mathbb{R}^d$ be a compact set. The following conditions are equivalent:

{\em (i)} $K$ is uniformly porous on lines.

{\em (ii)} No weak tangent of $K$ contains a line segment.

{\em (iii)} There exist constants $C > 0$ and $0 < \rho < 1$ such that 
$$
\mathcal N_{K \cap B(x,R) \cap T(L, r)}(r) \le C \biggl ( \frac{R}{r} \biggr )^{1-\rho},
$$
for every $x \in \mathbb{R}^d$, $0 < r < R$ and line $L$.
\end{lemma}

\begin{proof}[Proof of Theorem \ref{energy-comparison}]
(1) $\Rightarrow$ (2). Suppose $F \subset \gamma$, where $\gamma$ is a rectifiable curve. For a cube $Q \in \mathcal D$, let $3Q$ denote the cube with the same center and three times the side length. The {\em beta number} of $\gamma$ in $3Q$ is defined as
$$
\beta_\gamma(3Q) := \frac{1}{\diam(3Q)} \, \inf_L \sup_{z\in\gamma\cap 3Q}\dist(z,L),
$$
where the infimum is taken over all lines $L\subset\mathbb R^d$. For a constant $c > 0$ to be chosen, we call a cube $Q \in \mathcal D$ {\em wiggly} if $ \beta_\gamma(3Q) > c $ and {\em flat} otherwise. Define
$$
A_n := \sum_{\substack{Q \in \mathcal D_n \\ Q \text{ active}}} \diam Q
\qquad \text{and} \qquad
W_n := \sum_{\substack{Q \in \mathcal D_n \\ Q \text{ active and wiggly}}} \diam Q.
$$
By the traveling salesman theorem of P.~Jones,
$$
\sum_{Q \in \mathcal D, \, Q \cap \gamma \ne \emptyset}
\beta_\gamma(3Q)^2\diam Q
\lesssim \ell(\gamma),
$$
see \cite[Theorem 10.2.3]{bishop-peres} for the planar case and \cite{okikiolu} for the general case. Consequently,
$$
\sum_{n=0}^\infty W_n \lesssim \ell(\gamma).
$$

\smallskip 

\noindent {\em Step 1.}
Suppose that $Q\in\mathcal D_n$ is a flat active cube. By the definition of beta numbers, if $3c \le 2^{-m}$, then there is a line $L$ such that
$$
F\cap Q \, \subset \, \gamma \cap 3Q \, \subset \, T(L, 3c \cdot \diam Q)  \, \subset \, T(L, 2^{-m} \cdot \diam Q).
$$ 
Each descendant $Q' \subset Q$ in $\mathcal D_{m+n}$ has diameter $2^{-m} \diam Q$. Let $x$ be the center of $Q$. Applying the tube condition with $x, L$ and radii $R = \diam Q$, $r = 2^{-m} \diam Q$ shows that the number of these descendants which meet $F$ is at most
$$
C' \cdot 2^{m(1-\rho)}.
$$
Hence, the sum of their diameters is at most
$$
C' \cdot 2^{-m\rho}\diam Q.
$$
We first choose $m$ sufficiently large that this is at most $\frac12\diam Q$ and then choose $c>0$ so that $3c\le2^{-m}$. Consequently,
\begin{equation}
\label{eq:flat-active-cube-bound}
\sum_{\substack{Q' \subset Q, \, Q' \in \mathcal D_{m+n}, \\ Q' \text{ active} }} \diam Q' \le \frac{1}{2} \, \diam Q.
\end{equation}

For a wiggly active cube $Q \in \mathcal D_n$, we use the crude bound
\begin{equation}
\label{eq:wiggly-active-cube-bound}
\sum_{\substack{Q' \subset Q, \, Q' \in \mathcal D_{m+n}, \\ Q' \text{ active} }} \diam Q' \le M \diam Q,
\qquad M := 2^{(d-1)m}.
\end{equation}

\medskip

\noindent {\em Step 2.} 
Summing over the active cubes of generation $n$, we obtain
$$
A_{n+m}
\, \le \, \frac{1}{2} \cdot (A_n - W_n) + M W_n
\, \le \, \frac{1}{2} \cdot A_n+M W_n.
$$
Multiplying by 2 and summing over $0\le n \le N$, we get
$$
\sum_{n=0}^N A_n \, \le \, 2 \sum_{j=0}^{m-1} A_j + 2M \sum_{n=0}^\infty W_n  \, \lesssim \, 1 + \ell(\gamma).
$$
Letting $N \to \infty$, we obtain
\begin{equation}
\label{eq:rect-energy-estimate}
\mathscr E(F) \, = \, \sum_{n=0}^\infty A_n \, \lesssim \,1 + \ell(\gamma),
\end{equation}
where the implicit constant depends only on the dimension and the UPL constant of $K$.

\medskip

\noindent (2) $\Rightarrow$ (3) 
For $n \ge 0$, let $\mathcal S_n$ denote the collection of active cubes of side length $2^{-n}$. Since $\mathcal N(r) \lesssim \# \mathcal S_{n}$ for $2^{-n-1} < r \le 2^{-n}$,
$$
\int_0^1 \mathcal N(r) dr \, = \, \sum_{n=0}^\infty  \int_{2^{-n-1}}^{2^{-n}} \mathcal N(r) dr 
\, \lesssim \, \sum_{n=0}^\infty  2^{-n} \, \# \mathcal S_{n} \, \asymp \, \mathscr E(F).
$$

\smallskip

\noindent (3) $\Rightarrow$ (1) is immediate from Lemma \ref{jones-easy}.
\end{proof}

\section{Intersections with chord-arc curves}

In this section, we give a similar characterization of the closed subsets of a UPL set that are contained in chord-arc curves:

\begin{theorem}
\label{cac-energy-comparison}
Let $K$ be a compact UPL set in $[0,1]^d$, $d \ge 2$. For a closed set $F \subset K$, the following conditions are equivalent:

{\em (1)} $F$ is contained in a chord-arc curve.

{\em (2)}  There is a constant $C > 1$ such that for any $Q_0 \in \mathcal D$,
$$\mathscr E(F, Q_0) \, := \,  \sum_{\substack{Q \in \mathcal D, \, Q \subset Q_0, \\ Q \cap F \ne \varnothing}} \diam Q \, < \, C \diam Q_0.$$

{\em (3)} $\Adim F < 1$.
\end{theorem}

Together with Theorem \ref{continuum-frostman2}, the implication (1) $\Rightarrow$ (3) shows Theorem \ref{cac-sharp-gauges}.

\begin{proof}
(1) $\Rightarrow$ (2). Suppose $F$ is contained in a chord-arc curve $\gamma$. If $F \cap Q_0 = \varnothing$, then the claim is immediate. Otherwise, let $\gamma_0 \subset \gamma$ be the subarc of $\gamma$ from the first to the last point of $F \cap Q_0$. By the chord-arc property, $\ell(\gamma_0) \lesssim \diam Q_0$. After rescaling $Q_0$ to the unit cube, the estimate (\ref{eq:rect-energy-estimate}) from the previous section gives
$$
\mathscr E(F, Q_0) \, \lesssim \, \ell(\gamma_0) + \diam Q_0 \, \lesssim \, \diam Q_0.
$$

\medskip

\noindent (2) $\Rightarrow$ (3). Increasing $C$ if necessary, we may assume that $C > 2$. Fix an active cube $Q_0 \in\mathcal D_m$ and for each $k\ge0$, define
$$
A_k :=
\sum_{\substack{ Q \in \mathcal D_{m+k},\, Q \subset Q_0 \\ Q \cap F\neq\varnothing}} \diam Q
\qquad \text{and} \qquad
S_k : =\sum_{j=k}^{\infty}A_j.
$$
By Condition (2), we have
$$
S_k
\, = \, \sum_{\substack{ Q\in\mathcal D_{m+k},\, Q\subset Q_0\\ Q\cap F\neq\varnothing}} \mathscr E(F, Q)
\, \le \, C A_k.
$$
The inequalities
$$
S_{k+1} \, = \, S_k - A_k \, \le \, (1 - 1/C) S_k
$$
and
$$
S_0 \, = \, \mathscr E(F, Q_0) \, \le \, C \diam Q_0
$$
imply that
\begin{equation}
\label{eq:energy-tail-decay}
S_k \le C \theta^k \diam Q_0, \qquad \theta \, = \, 1 - 1/C \, \in \, (1/2,1).
\end{equation}

Fix an $0<r< \diam Q_0$ and choose $k \ge 1$ so that
$$
2^{-k} \diam Q_0 \, \le \, r \, < \, 2^{-(k-1)} \diam Q_0.
$$ 
The collection $\mathcal S_k$ of active generation $k$ descendants of $Q_0$ covers $F \cap Q_0$. Each cube in $\mathcal S_k$ has diameter $2^{-k} \diam Q_0$ and the sum of their diameters is at most $S_k$. Thus,
$$
\#\mathcal S_k  \, \le \, \frac{2^k S_k}{\diam Q_0} 
  \,  \le \, C (2\theta)^{k} \, = \, C \, 2^{ s k},
$$
for some $0 < s < 1$. As $2^k \asymp (\diam Q_0)/r$,
$$
\#\mathcal S_k \lesssim \biggl (\frac{\diam Q_0}{r} \biggr)^{s}.
$$
Since every cube $Q \in \mathcal S_k$ is contained in a ball of radius $r$,
\begin{equation}
\mathcal N_{F\cap Q_0}(r) \lesssim \biggl ( \frac{\diam Q_0}{r} \biggr )^{s} .
\end{equation}

Finally, fix a point $x \in K$ and $0 < r < R< \diam K$. Since $B(x,R)$ is contained in a uniformly bounded number of cubes
$Q_j \in\mathcal D$ with $\diam Q_j \asymp R$,
$$
\mathcal N_{F\cap B(x,R)}(r)
\, \lesssim \, \sum_j \mathcal N_{F\cap Q_j}(r)
\, \lesssim \, \biggl (\frac {R}{r} \biggr)^{s}.
$$
Thus, by the definition of the Assouad dimension, $\Adim F \le s < 1$.

\medskip

\noindent (3) $\Rightarrow$ (1) is immediate from Lemma \ref{badger-vellis}.
\end{proof}

\subsection*{Acknowledgements} The authors wish to thank Ben Jaye for introducing them to UPL sets. T.V. was partially supported by an AMS--Simons travel grant.

\subsection*{Use of AI}
The proofs of Theorems 1.1 and 1.2 were obtained with help of Claude Fable 5. GPT-6.1 Sol was used to polish the language and presentation. The authors checked the proofs and streamlined the exposition.

\bibliographystyle{amsplain}

 \medskip

\noindent Oleg Ivrii: {\tt ivrii@tauex.tau.ac.il}

\smallskip
\noindent Wladimir Schreiber Institute of Mathematical Sciences, Tel Aviv University, Tel Aviv, Israel

\bigskip

\noindent Truong Vu: {\tt  truongvu@math.ucla.edu}

\smallskip
\noindent University of California, Los Angeles, California, USA
 
\end{document}